# HAMMERSLEY'S PROCESS WITH SOURCES AND SINKS

By Eric Cator and Piet Groeneboom

*Delft University of Technology*

We show that, for a stationary version of Hammersley's process, with Poisson "sources" on the positive $x$-axis, and Poisson "sinks" on the positive $y$-axis, an isolated second-class particle, located at the origin at time zero, moves asymptotically, with probability 1, along the characteristic of a conservation equation for Hammersley's process. This allows us to show that Hammersley's process without sinks or sources, as defined by Aldous and Diaconis [*Probab. Theory Related Fields* **10** (1995) 199–213] converges locally in distribution to a Poisson process, a result first proved in Aldous and Diaconis (1995) by using the ergodic decomposition theorem and a construction of Hammersley's process as a one-dimensional point process, developing as a function of (continuous) time on the whole real line. As a corollary we get the result that $EL(t,t)/t$ converges to 2, as $t \to \infty$, where $L(t,t)$ is the length of a longest North-East path from $(0,0)$ to $(t,t)$. The proofs of these facts need neither the ergodic decomposition theorem nor the subadditive ergodic theorem. We also prove a version of Burke's theorem for the stationary process with sources and sinks and briefly discuss the relation of these results with the theory of longest increasing subsequences of random permutations.

**1. Introduction.** Let $L_n$ be the length of a longest increasing subsequence of a random permutation of the numbers $1, \ldots, n$, for the uniform distribution on the set of permutations. As an example, consider the permutation $(5, 3, 6, 2, 8, 7, 1, 4, 9)$. Longest increasing subsequences are $(3, 6, 7, 9)$, $(3, 6, 8, 9)$, $(5, 6, 7, 9)$ and $(5, 6, 8, 9)$. In this example the length of a longest increasing subsequence is equal to 4.

In Hammersley (1972) a discrete-time interacting particle process was introduced, which has at the $n$th step a number of particles equal to the length of a longest increasing subsequence of a (uniform) random permutation of length $n$. This process is defined in the following way.

Start with zero particles. At each step, let, according to the uniform distribution on $[0,1]$, a random particle $U$ in $[0,1]$ appear; simultaneously, let the nearest particle (if any) to the right of $U$ disappear. Then, as shown in Hammersley (1972), the number of particles after $n$ steps is distributed as $L_n$. Hammersley (1972) uses this discrete-time interacting particle process to show that $EL_n/\sqrt{n}$ converges to a finite constant $c > 0$, which is also the limit in probability [and, as noticed later by H. Kesten in his discussion of Kingman (1973), the almost sure limit] of $L_n/\sqrt{n}$. To prove that $EL_n/\sqrt{n}$ converges to a finite constant $c > 0$ is the first part of "Ulam's problem," the second part being the determination of $c$.

Aldous and Diaconis (1995) introduce a continuous-time version of the interacting particle process in Hammersley (1972), letting new particles appear according to a Poisson process of rate 1, using the following rule:

EVOLUTION RULE. At times of a Poisson (rate $x$) process in time, a point $U$ is chosen uniformly on $[0, x]$, independent of the past, and the particle nearest to the right of $U$ is moved to $U$, with a new particle created at $U$ if no such particle exists in $[0, x]$.

For our purposes the following alternative description is most useful. Start with a Poisson point process of intensity 1 on $\mathbb{R}_+^2$. Now shift the interval $[0, x]$ vertically through (a realization of) this point process, and, each time a point is caught, shift to this point the previously caught point that is immediately to the right. Let $L(x, y)$ be the number of particles in the interval $[0, x]$ after shifting to height $y$. Then, by Poissonization of the length of the random permutation, we get

$$L_{\widetilde{N}_{x,y}} \stackrel{\mathcal{D}}{=} L(x, y),$$

where

$$\widetilde{N}_{x,y} = \#\{\text{points of Poisson point process in } [0, x] \times [0, y]\} \stackrel{\mathcal{D}}{=} \text{Poisson}(xy).$$

In an alternative interpretation, $L(x, y)$ is the maximal number of points on a North-East path from $(0, 0)$ to $(x, y)$ with vertices at the points of the Poisson point process in the interior of $\mathbb{R}_+^2$, where the length of a North-East path is defined as the number of vertices it has at the points of the Poisson point process in the interior of $\mathbb{R}_+^2$. The reason is that a longest North-East path from the origin to $(x, y)$ has to pick up a point from each space–time path crossing the rectangle $[0, x] \times [0, y]$. Aldous and Diaconis (1995) call the evolving point process $y \mapsto L(\cdot, y)$, $y \geq 0$, of newly caught and shifted points *Hammersley's interacting particle process*.

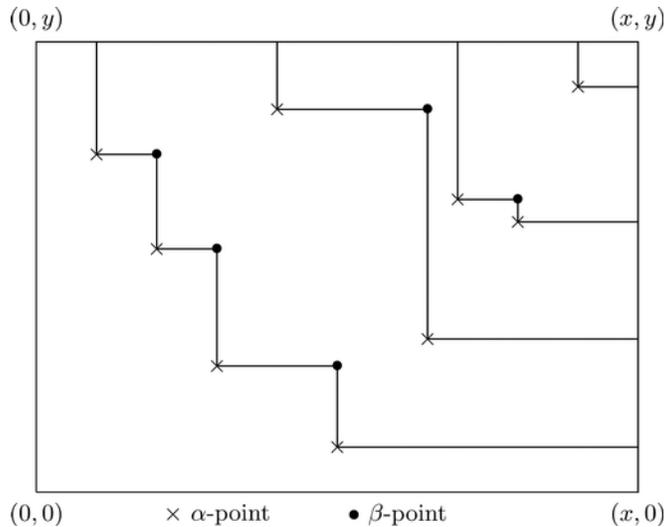

Fig. 1. *Space–time paths of Hammersley's process, contained in $[0,x] \times [0,y]$.*

We can also introduce the evolving point process $x \mapsto L(x,\cdot)$, $x \geq 0$, running from left to right. Analogously to the description above of the process running up, we shift in this case an interval $[0,y]$ on the $y$-axis to the right through the point process in the interior of the first quadrant, and, each time a point is caught, shift to this point the previously caught point that is immediately below this point (if there is such a point). By symmetry, it is clear that the processes $y \mapsto L(\cdot,y)$, $y \geq 0$, and $x \mapsto L(x,\cdot)$, $x \geq 0$, have the same distribution.

A picture of the space–time paths corresponding to the permutation $(5,3,6,2,8,7,1,4,9)$ is shown in Figure 1. In this case $[0,x] \times [0,y]$ contains nine points, and one can check graphically that there are four longest North-East paths (of length 4) from $(0,0)$ to $(x,y)$, corresponding to the subsequences $(3,6,7,9)$, $(3,6,8,9)$, $(5,6,7,9)$ and $(5,6,8,9)$. Following a terminology introduced in Groeneboom (2001), we call the points of the Poisson point process in the interior of $\mathbb{R}_+^2$ $\alpha$-*points* and the North-East corners of the space–time paths of Hammersley's process $\beta$-*points*. In fact, the actual $x$-coordinates of the $\alpha$-points in the picture are different from the numbers $3, 6, \ldots$, but the ranks of these $x$-coordinates are given by $3, 6$, and so on, if we order the $\alpha$-points according to the second coordinate.

We use a further extension of Hammersley's interacting particle process, where we have not only a Poisson point process in the interior of $\mathbb{R}_+^2$, but also, independently of this Poisson point process, mutually independent Poisson point processes on the $x$- and $y$-axis. We call the Poisson point process on the $x$-axis a process of "sources," and the Poisson point process on the $y$-axis a process of "sinks." The motivation for this terminology is that we now start the interacting particle process with a nonempty configuration of "sources" on the $x$-axis, which are subjected to the Hammersley's interacting particle process in the interior of $\mathbb{R}_+^2$, and which "escape" through sinks on the $y$-axis, if such a sink appears to the immediate left of a particle (with no other particles in between). Figure 2 shows how the space–time paths change if we add two sources and three sinks (at particular locations) to the configuration in Figure 1.

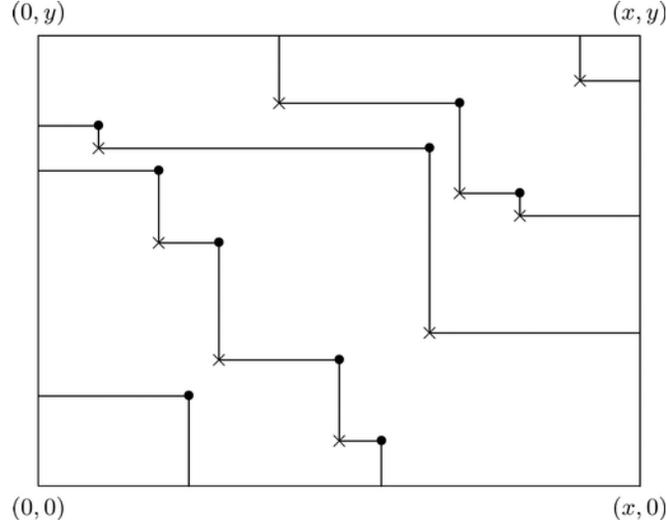

Fig. 2. *Space–time paths of Hammersley's process, with sources and sinks.*

The interacting particle process with sources and sinks was studied in Section 4 of Groeneboom (2002), where it was proved that, if the intensity of the Poisson processes on the $x$- and $y$-axes are $\lambda$ and $1/\lambda$, respectively, and the intensity of the Poisson process in the interior of $\mathbb{R}_+^2$ is 1, the process is stationary in the sense that the crossings of the space–time paths of the half-lines $\mathbb{R}_+ \times \{y\}$ are distributed as a Poisson point process of intensity $\lambda$, for all $y > 0$. The stationarity of the process was proved by an infinitesimal generator argument. It also follows from the computations in the Appendix of the present paper. The process is studied from an analytical point of view in Baik and Rains (2000) (see Remark 3.1 in Section 3).

In Section 2 we compare Hammersley's interacting particle process, as introduced in Aldous and Diaconis (1995), with the stationary extension of this process, with sources on the $x$-axis, and sinks on the $y$-axis. However, as an intermediate step, we introduce a process with Poisson sources on the positive $x$-axis, but no sinks on the $y$-axis. From Theorem 2.1 in the present paper we can deduce that this particle process, with Poisson sources of intensity $\lambda$ on the positive $x$-axis, but no sinks on the $y$-axis, behaves below an asymptotically linear "wave" of slope $\lambda^2$ through the $\beta$-points as a stationary process.

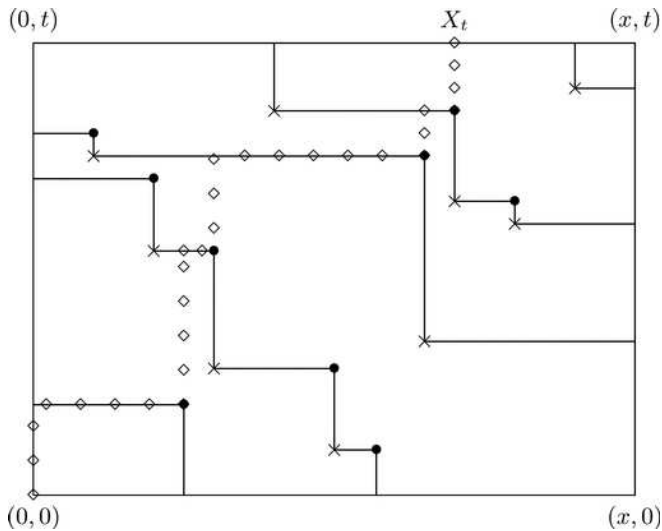

FIG. 3. *Path of isolated second-class particle in the configuration of Figure 2.*

In a coupling of the process with the stationary process, having both sources and sinks, this wave can be interpreted as the space–time path of an isolated second-class (or "ghost") particle with respect to the stationary process. For the concept "second-class particle" in the context of totally asymmetric simple exclusion processes (TASEP), see, for example, Ferrari (1992) or Liggett [(1999), Chapter 3]. The second-class particle jumps to the previous position of the particle that exits through the first sink at the time of exit, and successively jumps to the previous positions of particles directly to the right of it, at times where these particles jump to a position to the left of the second-class particle; see Figure 3. The space–time path of the isolated second-class particle moves asymptotically, with probability 1, along the characteristic of a conservation equation for the stationary process. Here we establish a connection with the theory of totally asymmetric simple exclusion processes. Although we use similar techniques as used for the study of the behavior of second-class particles in TASEP, the situation is in a certain sense simpler in our case, since we do not have to condition on having a second-class particle at the origin at time zero.

In a similar way we prove that Hammersley's process, with Poisson sinks of intensity $1/\lambda$, $\lambda > 0$, on the positive $y$-axis, but no sources on the $x$-axis, behaves asymptotically as a stationary process *above* a wave through the $\beta$-points of slope $\lambda^2$, if the Poisson sinks on the positive $y$-axis and the points of the Poisson process (of intensity 1) in the interior of $\mathbb{R}_+^2$ are independent. By a coupling argument, these processes can be compared directly to Hammersley's process, as defined in Aldous and Diaconis (1995), which has empty configurations on the $x$- and $y$-axis. The coupling argument gives a direct and "visual" proof of the local convergence of Hammersley's process to a Poisson point process with intensity $\lambda$, if one moves out along a "ray" $y = \lambda^2 x$, which is the main result Theorem 5 of Aldous and Diaconis (1995). The convergence of $EL(t,t)/t$ to 2, as $t \to \infty$, then also easily follows. This implies that $EL_n/\sqrt{n}$ converges to 2, a result first proved by Logan and Shepp (1977) and Vershik and Kerov (1977).

In Section 3 we study the $\beta$-points of the stationary Hammersley process. For these points we prove a "Burke theorem," showing that these points inherit the Poisson property from the $\alpha$-points. This allows us to show, using a time reversal argument, that in the stationary version of Hammersley's process, a longest "weakly" North-East path (allowing horizontal and vertical pieces along the $x$- or $y$-axis) only spends a vanishing fraction of time on the $x$- or $y$-axis.

**2. Path of an isolated second-class particle and local convergence of Hammersley's process.** Fix $\lambda > 0$, and let $t \mapsto L_\lambda(\cdot, t)$ be Hammersley's process, now considered as a one-dimensional point process, developing in time $t$, generated by a Poisson process of sources on the positive $x$-axis of intensity $\lambda$, $\lambda > 0$, a Poisson process of sinks on the time axis of intensity $1/\lambda$ and a Poisson process of intensity 1 in $\mathbb{R}_+^2$, where the Poisson process on the $x$-axis, the Poisson process on the time axis and the Poisson process in the plane are independent. It is helpful to switch from time to time the point of view of Hammersley's process as a process of space–time paths in $\mathbb{R}_+^2$ and Hammersley's process as a one-dimensional point process, developing in time. This is somewhat similar to the two ways one can view the Brownian sheet. Since the second coordinate can (mostly) be interpreted as "time" in the sequel, we will denote this coordinate by $t$ instead of $y$, although, with slight abuse of language, we will continue to call the vertical axis the "$y$-axis," following standard terminology.

We add an isolated second-class particle to the process, which is located at the origin at time zero. A picture of the trajectory of the isolated second-class particle for the configuration shown in Figure 2 is shown in Figure 3. Theorem 2.1 shows that the space–time path of the second-class particle is asymptotically linear with slope $\lambda^2$. This is to be expected from results on totally asymmetric simple exclusion processes (TASEP), as given in, for example, Ferrari (1992). For TASEP Burgers' equation is the relevant conservation equation in a continuous approximation to the process. The analogue of Burgers' equation for a macroscopic approximation to Hammersley's process (with neither sources nor sinks) is

$$(2.1) \qquad \frac{\partial u(x,t)}{\partial t} + u(x,t)^{-2}\frac{\partial u(x,t)}{\partial x} = 0,$$

where $u(x,t)$ is the intensity of the crossings at $(x,t)$; see Liggett [(1999), page 316], where the corresponding equation is given for the integrated intensity.

This leads us to expect that, analogously to the TASEP results,

$$t^{-1}X_t \xrightarrow{\text{a.s.}} 1/\lambda^2, \qquad t \to \infty,$$

where $X_t$ is the $x$-coordinate of the second-class particle, and where $\xrightarrow{\text{a.s.}}$ denotes almost sure convergence, since in this case the path $\{(x,t) = (t/\lambda^2, t) : t \geq 0\}$ is a characteristic for (2.1); compare to, for example, (12.1) in Section 12 of Ferrari (1992).

THEOREM 2.1. *Let $t \mapsto L_\lambda(\cdot, t)$ be the stationary Hammersley process, defined above, with intensities $\lambda$ and $1/\lambda$ on the $x$- and $y$-axis, respectively. Let $X_t$ be the $x$-coordinate of an isolated second-class particle w.r.t. $L_\lambda$ at time $t$, located at the origin at time zero. Then*

$$(2.2) \qquad t^{-1}X_t \xrightarrow{\text{a.s.}} 1/\lambda^2, \qquad t \to \infty.$$

The proof of Theorem 2.1 is based on Lemma 2.1. To formulate this lemma we first introduce some notation. Let $\eta_t, t \geq 0$, be the stationary point process, obtained by starting with a Poisson point process with intensity $\gamma > 0$ in $(0, \infty)$ at time 0, and letting it develop according to Hammersley's process on $(0, \infty)$, with Poisson sinks of intensity $1/\gamma$ on the $y$-axis, and a Poisson point process of intensity 1 in the interior of the first quadrant. Furthermore, let $\sigma_t, t \geq 0$, be the stationary process, coupled to $\eta_t, t \geq 0$, by using the same points in the first quadrant as used for $\eta$, and starting with a $(\delta/\gamma)$-"thickening," $\delta > \gamma$, of the Poisson point process with intensity $\gamma > 0$ on the $x$-axis, obtained by adding independently a Poisson point process of intensity $\delta - \gamma$, and letting $\sigma_t$ develop according to Hammersley's process on $(0, \infty)$. To get stationarity for the process $\sigma$, we replace the sinks on the $y$-axis by a $\gamma/\delta$-thinned set, obtained by keeping each sink with probability $\gamma/\delta$, independently for each sink. Then the sinks on the $y$-axis for the process $\sigma$ have intensity $1/\delta$. Finally, we let $t \mapsto \xi_t$ be the process of second-class particles of $\eta$ w.r.t. $\sigma$, that is, the points of $\xi_t$ denote the locations where the point process $\sigma_t$ has extra particles w.r.t. the point process $\eta_t$.

We use the notation $\eta_t[0, x]$ for the number of particles of $\eta_t$ in the interval $[0, x]$ at time $t$, with the convention that particles, escaping through a sink in the time interval $[0, t]$, are located at zero. We define $\sigma_t[0, x]$ similarly. Furthermore, we use the notation $\eta_t(0, x]$ ($\sigma_t(0, x]$) for the number of particles of $\eta_t$ ($\sigma_t$) in the open half-open interval $(0, x]$ at time $t$. Finally we define the "flux" $F_\xi(x, t)$ of $\xi$ through $x$ at time $t$ by

$$(2.3) \qquad F_\xi(x, t) = \sigma_t[0, x] - \eta_t[0, x].$$

FIG. 4.  *Processes $\eta$ and $\xi$.*

The flux $F_\xi(x,t)$ is equal to the number of second-class particles in $(0,x]$ at time $t$ minus the number of removed sinks in the segment $\{0\} \times [0,t]$ (through which space–time paths of second-class particles start moving to the right). Relation (2.3) is in fact a conservation law.

A picture of the processes $\eta$ and $\xi$ is shown in Figure 4. In this case the process $\sigma$ (inside the rectangle $[0,x] \times [0,t]$) is obtained from the process $\eta$ by adding two sources at the locations $z_1(0)$ and $z_2(0)$ and removing a sink at height $S_0$. The crossings of horizontal lines of the space–time paths of the process $\sigma$ are the unions of the crossings of (the same) horizontal lines of the space–time paths of the processes $\eta$ and $\xi$.

LEMMA 2.1. (i) *Let $\eta$ be Hammersley's process, defined above, with sources of intensity $\gamma > 0$ and sinks of intensity $1/\gamma$, and let $\delta > \gamma$. We add independently a Poisson point process of intensity $\delta - \gamma$ to the Poisson process of sources, and perform a $\gamma/\delta$-thinning of the Poisson point process of sinks of intensity $1/\gamma$ on the y-axis. Let $\sigma$ be Hammersley's process, coupled to $\eta$, and having the augmented set of sources with intensity $\delta$ and the thinned set of sinks with intensity $1/\delta$. Finally, let $Z_t$ be, at time $t$, the location of the second-class particle for which the space–time path starts moving to the right through the smallest removed sink. Then*

$$\lim_{t \to \infty} \frac{Z_t}{t} = \frac{1}{\gamma \delta} \qquad a.s.$$

(ii) *Let $\eta'$ represent Hammersley's process developing from left to right, with sources (on the x-axis) of intensity $\gamma > 0$ and sinks (on the y-axis) of intensity $1/\gamma$, and let $0 < \delta < \gamma$. We add independently a Poisson point process of intensity $\delta^{-1} - \gamma^{-1}$ to the Poisson process of sinks of intensity $\gamma^{-1}$, and perform a $\delta/\gamma$-thinning of the Poisson point process of sources of intensity $\gamma$ on the x-axis. Let $\sigma'$ be the process developing from left to right, coupled to $\eta'$, and having the augmented set of sinks with intensity $\delta^{-1}$ as sources and the thinned set of sources with intensity $\delta$ as sinks. Finally, let $Z'_t$ be the location of the second-class particle of $\sigma'$ w.r.t. $\eta'$, for which the space–time path leaves the x-axis through the smallest removed source (of the original process $\eta$). Note that the smallest removed source of $\eta$ is a removed sink for $\eta'$. Then*

$$\lim_{t \to \infty} \frac{Z'_t}{t} = \gamma\delta \qquad a.s.$$

PROOF. (i) Let $x > 0$. We have

$$\lim_{n \to \infty} \frac{\eta_n[0, nx]}{n} = \frac{1}{\gamma} + x\gamma \qquad \text{a.s.},$$

since $\eta_n[0, nx]$ equals $\eta_n(0, nx]$ plus the number of sinks for the process $\eta$, contained in $\{0\} \times [0, n]$ (where $n$ is a positive integer), and since $\eta_n(0, nx]$ and the number of sinks contained in $\{0\} \times [0, n]$ have Poisson distributions with parameters $nx\gamma$ and $n/\gamma$, respectively. Here we use the stationarity of the process $\eta$, implying that $\eta_n(0, nx]$ has a Poisson distribution with parameter $nx\gamma$. Note that, for each $\varepsilon > 0$,

$$\sum_{n=1}^{\infty} P\{|\eta_n(0, nx] - nx\gamma| > n\varepsilon\} < \infty,$$

and hence, by the Borel–Cantelli lemma,

$$P\{|\eta_n(0, nx] - nx\gamma| > n\varepsilon \text{ infinitely often}\} = 0,$$

implying the almost sure convergence of $\eta_n(0, nx]/n$ to $x\gamma$, as $n \to \infty$. The almost sure convergence to $1/\gamma$ of the number of sinks for the process $\eta$, contained in $\{0\} \times [0, n]$, divided by $n$, follows in the same way.

Similarly,

$$\lim_{n \to \infty} \frac{\sigma_n[0, nx]}{n} = \frac{1}{\delta} + x\delta \qquad \text{a.s.}$$

Hence, by (2.3),

$$(2.4) \qquad \lim_{n \to \infty} \frac{F_\xi(nx, n)}{n} = \frac{1}{\delta} - \frac{1}{\gamma} + x(\delta - \gamma) = -(\delta - \gamma)\left\{\frac{1}{\gamma\delta} - x\right\} \qquad \text{a.s.}$$

This limit is negative for $0 < x < 1/(\gamma\delta)$ and positive for $x > 1/(\gamma\delta)$.

We can number the particles of $\xi$ according to their position at time 0, so that, for $i > 0$, particle $i$ is the $i$th second-class particle to the right of the origin at time 0. We then let $z_i(t)$ be the position of the $i$th second-class particle at time $t \geq 0$. For $i \leq 0$, we let $z_i(t)$, $i = 0, -1, -2, \ldots$, be the second-class particles at time $t$, for which the space–time paths leave the $y$-axis through the removed sinks $S_0, S_1, \ldots$, respectively, ordering these removed sinks according to the height of their location on the $y$-axis; note that $Z_t = z_0(t)$ (see Figure 4).

Hence $F_\xi(x, t)$ has the representation

$$(2.5) \qquad F_\xi(x,t) = \#\{i > 0 : z_i(t) \leq x\} - \#\{i \leq 0 : z_i(t) > x\}.$$

Note that second-class particles $z_i(\cdot)$, $i \leq 0$, starting their space–time path to the right at a removed source in $\{0\} \times [0, t]$, and satisfying $z_i(t) \in [0, x]$, do not give a contribution to (2.5), since they give a contribution to $\eta_t[0, x]$ as a particle of $\eta_t$, located at zero, and a contribution to $\sigma_t[0, x]$ as a particle of $\sigma_t$ in the interval $(0, x]$. These two contributions cancel in (2.3). It is also clear from (2.5) that, for fixed $t$, the flux $F_\xi(x, t)$ is nondecreasing in $x$.

Relation (2.5) shows that $F_\xi(Z_n, n) = F_\xi(z_0(n), n)$ is equal to zero at each time $n$, and since $F_\xi(nx, n)$ is nondecreasing in $x$ for fixed $n$, we get from (2.4),

$$\lim_{n \to \infty} \frac{Z_n}{n} = \frac{1}{\gamma\delta} \qquad \text{a.s.}$$

But, since $Z_t$ is nondecreasing in $t$, we then also have

$$\lim_{t \to \infty} \frac{Z_t}{t} = \frac{1}{\gamma\delta} \qquad \text{a.s.}$$

(ii) The result is obtained from part (i) by reflecting the processes w.r.t. the diagonal, and noting that the reflected processes have the same probabilistic behavior, but with the role of sources and sinks interchanged. The limit $1/(\gamma\delta)$ changes to $\gamma\delta$ because of the interchange of $x$- and $y$-coordinate. □

PROOF OF THEOREM 2.1. We couple the process $t \mapsto (L_\lambda(\cdot, t), X_t)$ with the process $t \mapsto (\eta_t, \sigma_t)$, where the processes $\eta$ and $\sigma$ are defined as in part (i) of Lemma 2.1, and where $L_\lambda(\cdot, t) = \eta_t$ and $\delta > \gamma = \lambda$. Then $Z_t \leq X_t$, for all $t \geq 0$, where $Z_t$ is defined as in part (i) of Lemma 2.1. This is seen in the following way.

At time zero, we have $Z_0 = X_0 = 0$. Since the process $\sigma$ is obtained from the process $\eta$ by a thinning of the sinks and a "thickening" of the sources, and the space–time path of $Z_t$ leaves the axis $\{0\} \times \mathbb{R}_+$ through the smallest removed sink, it will leave this axis at a time which is larger than or equal to the time the space–time path of $X_t$ leaves the axis, since the space–time path of $X_t$ will leave the axis through the smallest sink in the original set of sinks. Note that since $\sigma$ has less sinks and more sources:

(2.6) $$\eta_t(0, x] \leq \sigma_t(0, x], \qquad t \geq 0, x > 0.$$

This means that not only $Z_t$ becomes positive at a time that is at least as large as the time that $X_t$ becomes positive, but also moves to the right at a speed that is not faster than that of $X_t$. Also note that if $Z_t$ jumps to a position $x > Z_{t-}$, an $\eta$-particle jumps over it from a position $x' \geq x$. Here and in the sequel we use the notation $Z_{t-}$ to denote $\lim_{t' \uparrow t} Z_{t'}$, with a similar convention for $X_{t-}$.

If $X_{t-} < x$ and $Z_{t-} \leq X_{t-}$, $X_t$ will jump to $x'$. Since $Z_t \leq x'$, $Z_t$ can never overtake $X_t$. Note that we can have $x' > x$ if several second-class particles are next to each other, without a first-class particle in between. In this case $Z_t$ does not have to move to the position of the $\eta$ particle, but can move to the position of the closest second-class particle to the right of it.

Hence we have, with probability 1,

$$\liminf_{t \to \infty} \frac{X_t}{t} \geq \lim_{t \to \infty} \frac{Z_t}{t} = \frac{1}{\gamma \delta} = \frac{1}{\delta \lambda}.$$

Since this is true for any $\delta > \lambda$, we get

$$\liminf_{t \to \infty} \frac{X_t}{t} \geq \frac{1}{\lambda^2}.$$

For the reverse inequality, we switch the role of the sources and the sinks, and view Hammersley's process as developing from left to right. This time we add independently a Poisson point process of intensity $\delta^{-1} - \gamma^{-1}$ to the Poisson process of sinks of intensity $\gamma^{-1}$, and perform a $\delta/\gamma$-thinning of the Poisson point process of sources of intensity $\gamma$ on the $x$-axis, where $\gamma = \lambda$ and $0 < \delta < \gamma$, and use the process $\eta'$ and $\sigma'$, defined in part (ii) of Lemma 2.1. Note that $\eta'$ has the same space–time paths as the process $\eta$, defined above. In the coupling we now consider $L_\lambda$ as a process developing from left to right and take $L_\lambda(t, \cdot) = \eta'_t$.

Let $X'_x$ be an isolated second-class particle for the process running from left to right in the same way as $X_t$ is an isolated second-class particle for the process running upward. Trajectories of $X$ and $X'$ are shown in Figure 5.

We have

(2.7) $$X(X'(x)) \leq x, \qquad x \geq 0,$$

writing temporarily $X'(x)$ instead of $X'_x$ and $X(u)$ instead of $X_u$. Equation (2.7) is equivalent to noting that the trajectory of $(X_t, t)$ lies above the trajectory of $(x, X'_x)$ (see also Figure 5). This follows from the fact that if $(X_t, t)$ hits a space–time path at a point North-West of the point where $(x, X'_x)$ hits the same space–time path, this must also be true for the next space–time path, since the first trajectory moves up, and the second trajectory moves to the right.

By Lemma 2.1 and the argument above, now applied on the process moving from left to right, we get the relation

$$\liminf_{x \to \infty} \frac{X'_x}{x} \geq \lim_{x \to \infty} \frac{Z'_x}{x} = \delta\lambda, \tag{2.8}$$

with probability 1. But the almost sure relation $\liminf_{x \to \infty} X'_x/x \geq \delta\lambda$ implies for the process $t \mapsto X_t$ the almost sure relation

$$\limsup_{t \to \infty} \frac{X_t}{t} \leq 1/(\delta\lambda), \tag{2.9}$$

since we get for each $\lambda' > 1/(\delta\lambda)$, with probability 1,

$$\limsup_{t \to \infty} \frac{X(t/\lambda')}{t/\lambda'} \leq \limsup_{t \to \infty} \frac{X(X'(t))}{t/\lambda'} \leq \lim_{t \to \infty} \frac{t}{t/\lambda'} = \lambda',$$

using (2.8) in the first inequality and (2.7) in the second inequality.

Since (2.9) is true for any $\delta < \lambda$, we get, with probability 1,

$$\limsup_{t \to \infty} \frac{X_t}{t} \leq \frac{1}{\lambda^2}.$$

The result now follows. □

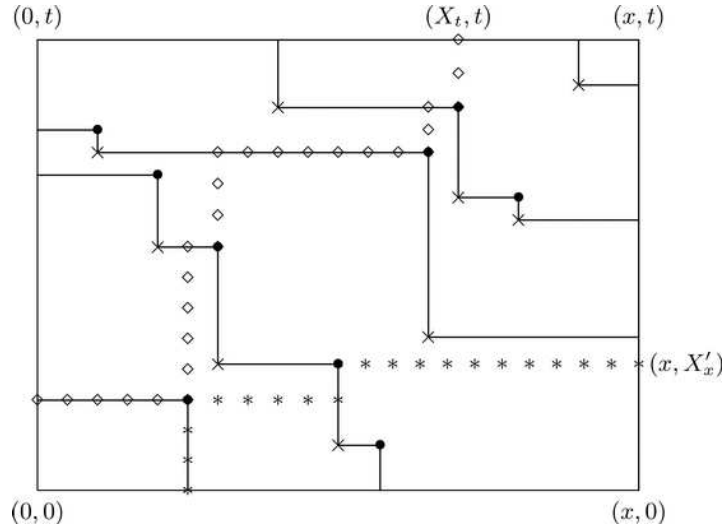

FIG. 5. *Trajectories of $(X_t, t)$ and $(x, X'_x)$.*

REMARK 2.1. The second-class particle $X'_x$, introduced at the end of the proof of Theorem 2.1, plays the same role for Hammersley's process, running from left to right, as the second-class particle $X_t$ plays for Hammersley's process, running up. It therefore has to satisfy

$$\lim_{x \to \infty} \frac{X'_x}{x} = \lambda^2, \tag{2.10}$$

with probability 1. Note that we get an interchange of the $x$ and $t$ coordinate which leads to $\lambda^2$ in (2.10) instead of the $1/\lambda^2$ in (2.2), but that the line along which $(x, X'_x)$ tends to $\infty$ is in fact the same as the line along which $(X_t, t)$ tends to $\infty$.

The following lemma will allow us to show that Theorem 2.1 implies both the local convergence of Hammersley's process to a Poisson process and the relation $c = 2$ [which is the central result Theorem 5 on page 204 in Aldous and Diaconis (1995)].

LEMMA 2.2. *Let $L_\lambda$ be the stationary Hammersley process, defined in Theorem 2.1. Furthermore, let $L_\lambda^{-y}$ be the process obtained from $L_\lambda$ by omitting the sinks on the $y$-axis, and let $L_\lambda^{-x}$ be the process obtained from $L_\lambda$ by omitting the sources on the $x$-axis. $L_\lambda^{-y}$ is coupled to $L_\lambda$, by using the same point process in the interior of $\mathbb{R}_+^2$, and the same set of sources on the $x$-axis, and $L_\lambda^{-x}$ is coupled to $L_\lambda$, by using the same point process in the interior of $\mathbb{R}_+^2$, and the same set of sinks on the $y$-axis. Then:*

(i) *The processes $L_\lambda$ and $L_\lambda^{-y}$ have the same space–time paths below the space–time path $t \mapsto (X_t, t)$ of the isolated second-class particle $X_t$ for the process $t \mapsto L_\lambda(\cdot, t)$.*

(ii) *The processes $L_\lambda$ and $L_\lambda^{-x}$ have the same space–time paths above the space–time path $t \mapsto (t, X'_t)$ of the isolated second-class particle $X'_t$ for the process $t \mapsto L_\lambda(t, \cdot)$, running from left to right.*

PROOF. Omit the first sink at location $y_1$ on the $y$-axis. Then the path of $L_\lambda$ leaving through $(0, y_1)$ is changed to a path traveling up through the $\beta$-point with $y$-coordinate $y_1$ to the right of $(0, y_1)$ until it hits the next path of the original process. At this level the path of the changed (by omitting the smallest sink) process is going to travel to the left, and the next path will go up (instead of to the left) through the closest $\beta$-point to the right. And so on. The "wave" through the $\beta$-points that is caused by leaving out the first sink is in fact the space–time path of the isolated second-class particle $X_t$ (see Figure 3).

We can now repeat the argument for the situation that arises by leaving out the second sink. This will lead to a "wave" through $\beta$-points that is going to travel North of the first wave that was caused by leaving out the first sink. This wave is the space–time path of an isolated second-class particle in the new situation, where the first sink is removed. Below the first wave the space–time paths remain unchanged. The argument runs the same for all the remaining sinks.

(ii) The argument is completely similar, but now applies to the process running from left to right instead of up (see the end of the proof of Theorem 2.1).  □

In the proof of Corollary 2.1 we will need the concept of a "weakly North-East path," a concept also used in Baik and Rains (2000).

DEFINITION 2.1. In the stationary version of Hammersley's process, a *weakly North-East path* is a North-East path that is allowed to pick up points from either the Poisson process on the $x$-axis or the Poisson process on the $y$-axis before going strictly North-East, picking up points from the Poisson point process in the interior $\mathbb{R}_+^2$. The *length of a weakly North-East path* from $(0,0)$ to $(x,t)$ is the number of points of the Poisson processes on the axes and the interior of $\mathbb{R}_+^2$ on this path from $(0,0)$ and $(x,t)$. A *strictly North-East path* is a path that has no vertical or horizontal pieces (and hence no points from the axes).

Note that the length of a longest weakly North-East path from $(0,0)$ to $(x,t)$ in the stationary version of Hammersley's process is equal to the number of space–time paths intersecting $[0,x] \times [0,t]$, just as in the case of Hammersley's process without sources or sinks (in which case only strictly North-East paths are possible).

COROLLARY 2.1 [Theorem 5 of Aldous and Diaconis (1995)]. *Let $L$ be Hammersley's process on $\mathbb{R}_+$, started from the empty configuration on the axes. Then:*

(i) *For each fixed $a > 0$, the random particle configuration with counting process*

$$y \mapsto L(t+y, at) - L(t, at), \qquad y \geq -t,$$

*converges in distribution, as $t \to \infty$, to a homogeneous Poisson process on $\mathbb{R}$, with intensity $\sqrt{a}$.*

(ii)

$$\lim_{t \to \infty} EL(t,t)/t = 2.$$

PROOF. (i) Fix $a' > a$, and let, for $\lambda = \sqrt{a'}$, $L_\lambda^{-y}$ be Hammersley's process, starting from Poisson sources of intensity $\lambda$ on the positive $x$-axis, and running through an independent Poisson process of intensity 1 in the plane (without sinks). Then we get from Theorem 2.1 and Lemma 2.2 that the counting process $y \mapsto L_\lambda^{-y}(t+y, at) - L_\lambda^{-y}(t, at)$ converges in distribution to a Poisson process of intensity $\lambda$, since the process, restricted to a finite interval, lies with probability 1 at level $t$ to the right of the space–time path of the isolated second-class particle $X_t$, as $t \to \infty$.

If we couple the original Hammersley process and the process $L_\lambda^{-y}$ via the same Poisson point process in the plane, we get that at any level the number of crossings of horizontal lines of the process $L$ is contained in the set of crossings of these lines of the process $L_\lambda^{-y}$, since the latter process has sources on the $x$-axis and no sinks on the $y$-axis. Hence, for a finite collection of disjoint intervals $[a_i, b_i)$, $i = 1, \ldots, k$, and nonnegative numbers $\theta_1, \ldots, \theta_k$, we obtain

$$E \exp\left\{-\sum_{i=1}^k \theta_i \{L(t+b_i, at) - L(t+a_i, at)\}\right\}$$
$$\geq E \exp\left\{-\sum_{i=1}^k \theta_i \{L_\lambda^{-y}(t+b_i, at) - L_\lambda^{-y}(t+a_i, at)\}\right\}.$$

But the right-hand side converges by Theorem 2.1 and Lemma 2.2 to

$$\exp\left\{-\sum_{i=1}^k \lambda(b_i - a_i)\{1 - e^{-\theta_i}\}\right\},$$

so we get

(2.11)
$$\liminf_{t \to \infty} E \exp\left\{-\sum_{i=1}^k \theta_i \{L(t+b_i, at) - L(t+a_i, at)\}\right\}$$
$$\geq e^{-\sum_{i=1}^k \lambda(b_i - a_i)\{1 - e^{-\theta_i}\}}.$$

A similar argument, but now comparing the process $L$ with a process $L_\lambda^{-x}$, having sinks of intensity $1/\lambda = 1/\sqrt{a'}$ on the $y$-axis (which can be considered to be "sources" for Hammersley's process, running from left to right), but no sources on the $x$-axis, shows

(2.12)
$$\limsup_{t \to \infty} E \exp\left\{-\sum_{i=1}^k \theta_i \{L(t+b_i, at) - L(t+a_i, at)\}\right\}$$
$$\leq e^{-\sum_{i=1}^k \lambda(b_i - a_i)\{1 - e^{-\theta_i}\}},$$

for any $a' < a$, since in this case the crossings of horizontal lines of the process $L$ are supersets of the crossings of these lines by the process $L_\lambda^{-x}$.

That the crossings of horizontal lines of the process $L$ are supersets of the crossings of horizontal lines by the process $L_\lambda^{-x}$ can be seen in the following way. Proceeding as in the proof of Lemma 2.2, we can, for the process $L_\lambda$, omit the sources one by one, starting with the smallest source. The omission of the smallest source will generate the path of a second-class particle $X'_t$, and the paths of $L_\lambda$ will, at the interior of a vertical segment of the path of $X'_t$, have an extra crossing of horizontal lines w.r.t. the paths of the process with the omitted source. On the other hand, the process with the omitted source will have extra crossings of *vertical* lines, since some particles will make bigger jumps to the left. We can now repeat the argument by omitting the second source, which will lead to a further decrease of crossings of horizontal lines, and so on.

Combining (2.11) and (2.12), we find

$$\lim_{t\to\infty} E\exp\left\{-\sum_{i=1}^{k}\theta_i\{L(t+b_i,at)-L(t+a_i,at)\}\right\} = e^{-\sum_{i=1}^{k}(b_i-a_i)\sqrt{a}\{1-e^{-\theta_i}\}},$$

and the result follows.

(ii) Since the length of a longest strictly North-East path is always smaller than or equal to the length of a longest weakly North-East path, in the situation of a stationary process with Poisson sources on the positive $x$-axis and Poisson sinks on the positive $y$-axis, both with intensity 1, we must have, for each $t > 0$,

$$EL(t,t)/t \leq 2,$$

since the expected length of a longest weakly North-East path from $(0,0)$ to $(t,t)$ is $2t$ for the stationary process.

The latter fact was proved in Groeneboom (2002), and comes from the simple observation that the length of a longest weakly North-East path from $(0,0)$ to $(t,t)$ is equal to the total number of paths crossing $\{0\} \times [0,t]$ and $[0,t] \times \{t\}$. Since the number of crossings of $\{0\} \times [0,t]$ has a Poisson($t$) distribution by construction, and the number of crossings of $[0,t] \times \{t\}$ also has a Poisson($t$) distribution, this time by the stationarity of the process $L_\lambda$, where $\lambda = 1$ in the present case, we get that the expectation of the total number of crossings of the left and upper edge is exactly $2t$.

To prove conversely that $\liminf_{t\to\infty} EL(t,t)/t \geq 2$, we first note that $L(t,t)$ is in fact the number of crossings of Hammersley's space–time paths with the line segment $[0,t] \times \{t\}$. Take a partition $0, t/k, 2t/k, \ldots, t$ of the interval $[0,t]$, for some integer $k > 0$. Then the crossings of the space–time paths of $L$ of the segment $[(i-1)t/k, it/k] \times \{t\}$ *contain* the crossings of this line segment by the paths of a Hammersley process $L_{\lambda_i}^{-x}$ with sinks of intensity $1/\lambda_i = 1/\sqrt{a_i}$, $a_i < k/i$, on the $y$-axis, but no sources on the $x$-axis.

But, by Theorem 2.1 and Lemma 2.2, the crossings of the process $L_{\lambda_i}^{-x}$ with the segment $[(i-1)t/k, it/k] \times \{t\}$ belong, as $t \to \infty$, to the stationary part of the process with probability 1, since $a_i < k/i$.

We now have

$$\lim_{t \to \infty} t^{-1} E\{L_{\lambda_i}^{-x}(it/k, t) - L_{\lambda_i}^{-x}((i-1)t/k, t)\} = \frac{\lambda_i}{k},$$

by uniform integrability of $t^{-1} L_{\lambda_i}^{-x}(\gamma t, t)$, $\gamma \in (0, i/k]$, $t \geq 0$, using, for example, the fact that the second moments are bounded above by the second moments of the corresponding stationary process with sources of intensity $\lambda_i$ and sinks of intensity $1/\lambda_i$. Hence we get, by summing over the intervals of the partition,

$$\liminf_{t \to \infty} EL(t,t)/t \geq \frac{1}{k} \sum_{i=1}^{k} \sqrt{a_i}.$$

Letting $a_i \uparrow k/i$, we obtain (still for fixed $k$)

$$\liminf_{t \to \infty} EL(t,t)/t \geq \sum_{i=1}^{k} 1/\sqrt{ik} = 2(1 + O(1/k)),$$

and (ii) follows by letting $k \to \infty$ in the latter relation. □

**3. Burke's theorem for Hammersley's process.** In this section we show that, in the stationary version of Hammersley's process with sources on the $x$-axis and sinks on the $y$-axis, the $\beta$-points inherit the Poisson property from the $\alpha$-points. One could consider this as a version of Burke's theorem for Hammersley's process. Burke's theorem [see Burke (1956)] states that the output of a stationary $M/M/1$ queue is Poisson. An interesting generalization of Burke's theorem is discussed in O'Connell and Yor (2002). A version of Burke's theorem for totally asymmetric simple exclusion processes is given in Ferrari [(1992), Theorem 7.1]. Burke's theorem is essentially based on a time-reversibility property and for our result on the $\beta$-points this is also the case. Our version of Burke's theorem runs as follows.

THEOREM 3.1. *Let $L_\lambda$ be a stationary Hammersley process on $[0, T_1] \times [0, T_2]$, generated by a Poisson process of "sources" of intensity $\lambda$ on the positive $x$-axis, a Poisson process of intensity $1/\lambda$ of "sinks" on the positive $y$-axis and a Poisson process of intensity $1$ in $\mathbb{R}_+^2$, where the three Poisson processes are independent. Let $L_\lambda^\beta$ denote the point process of $\beta$-points in $[0, T_1] \times [0, T_2]$, that is, the North-East corners of the space–time paths of the process $L_\lambda$, restricted to $[0, T_1] \times [0, T_2]$, $L_\lambda^{\text{in}}$ the entries of the space–time* paths on the East side of $[0, T_1] \times [0, T_2]$ and $L_\lambda^{\text{out}}$ the exits of the space–time paths on the North side. Then $L_\lambda^\beta$ is a homogeneous Poisson point process with intensity 1 in $[0, T_1] \times [0, T_2]$, $L_\lambda^{\text{in}}$ is a homogeneous Poisson process of intensity $1/\lambda$ and $L_\lambda^{\text{out}}$ is a homogeneous Poisson process of intensity $\lambda$, and all three processes are independent.

PROOF. We define a state space $E$ as the possible finite point configurations on $[0, T_1]$, so $E = \bigsqcup_{n=0}^\infty E_n$, where

$$E_n = \{(x_1, \ldots, x_n) : 0 \leq x_1 \leq \cdots \leq x_n \leq T_1\} \qquad (n \geq 1)$$

and $E_0 = \{\varnothing\}$, the empty configuration. We endow each $E_n$ with the usual topology, which makes $E$ into a locally compact space. We define a Markov process $(X_t)_{0 \leq t \leq T_2}$ on $E$ such that $X_t$ is the point configuration of the Hammersley process $L$ on the line $[0, T_1] \times \{t\}$. In particular we have that $X_0$ is distributed according to a Poisson process with intensity $\lambda$. From the definition of the Hammersley process it is not hard to see that the generator $G$ of this Markov process is given by

$$Gf(x) = \int_0^{T_1} f(\mathcal{R}_t x)\, dt + \frac{1}{\lambda} f(\mathcal{L}x) - \left(\frac{1}{\lambda} + T_1\right) f(x)$$

where $f \in C_0(E)$, $\mathcal{L}$ corresponds to an exit to the left and $\mathcal{R}_t$ corresponds to an insertion of a new Poisson point at $t$, so

$$\mathcal{L} : E \to E : \mathcal{L}x = \begin{cases} (x_2, \ldots, x_n), & \text{if } x \in E_n \ (n \geq 2), \\ \varnothing, & \text{if } x \in E_0 \sqcup E_1, \end{cases}$$

and for $0 < t < T_1$,

$$\mathcal{R}_t : E \to E : \mathcal{R}_t x = \begin{cases} (x_1, \ldots, x_{i-1}, t, x_{i+1}, \ldots, x_n), \\ \qquad \text{if } x_{i-1} < t \leq x_i \ (x \in E_n), \\ (x_1, \ldots, x_n, t), \quad \text{if } x_n < t \ (x \in E_n). \end{cases}$$

Here we use the convention that $x_0 = 0$. To prove that $G$ is indeed the generator, we fix $f \in C_0(E)$ and $x \in E$ and consider the transition operators

$$P_t f(x) = E(f(X_t) | X_0 = x) \qquad (t \geq 0).$$

We will consider the process for a time interval $[0, h]$ $(h \downarrow 0)$ and call $A_h$ the number of Poisson points in the strip $[0, T_1] \times [0, h]$ and $S_h$ the number of sinks in $\{0\} \times [0, h]$. Then

$$P_h f(x) = f(x) P(A_h = 0 \text{ and } S_h = 0)$$
$$+ \frac{1}{T_1} \int_0^{T_1} f(\mathcal{R}_t x)\, dt \cdot P(A_h = 1 \text{ and } S_h = 0)$$
$$+ f(\mathcal{L}x) P(A_h = 0 \text{ and } S_h = 1) + O(h^2)$$
$$= f(x)\left(1 - T_1 h - \frac{1}{\lambda} h\right) + h \int_0^{T_1} f(\mathcal{R}_t x)\, dt + \frac{h}{\lambda} f(\mathcal{L}x) + O(h^2).$$

This shows that for every $f \in C_0(E)$ and every $x \in E$,

$$\frac{d}{dt}\bigg|_{t=0} P_t f(x) = Gf(x).$$

Since $X_t$ is clearly a homogeneous Markov process, we get for $t \in [0, T_2]$,

(3.1) $$\frac{d}{ds}\bigg|_{s=t} P_s f(x) = GP_t f(x).$$

Now we note that $G$ is a continuous operator on $C_0(E)$, so $e^{tG}$ exists and is also a continuous operator. Since

$$\frac{d}{ds}\bigg|_{s=t} e^{sG} f(x) = Ge^{tG} f(x),$$

(3.1) together with the uniqueness of solutions of a differential equation proves that

$$P_t f(x) = e^{tG} f(x).$$

The key idea to prove the theorem is to consider the time-reversed process

$$\widetilde{X}_s = \lim_{s' \downarrow s} X_{T_2 - s'} \qquad (\widetilde{X}_{T_2} = X_0).$$

We take the left-limit of the original process $X$ to ensure the càdlàg property of $(\widetilde{X}_s)_{0 \leq s \leq T_2}$. Since, given $X_t$, the past of the process $X$ is independent of the future, it follows immediately that $\widetilde{X}$ is a Markov process, possibly inhomogeneous. However, if we define $\mu$ as the probability measure on $E$ induced by a Poisson process of intensity $\lambda$, then $X_0 \sim \mu$ and $\mu$ is a stationary measure for the generator $G$, which implies that $\widetilde{X}$ also is stationary and homogeneous. The stationarity of $X$ was shown in Groeneboom (2002), but will also be a consequence of calculations done in the Appendix. Now consider the transition operators

$$\widetilde{P}_t f(x) = E(f(\widetilde{X}_t) | \widetilde{X}_0 = x) \qquad (t \geq 0)$$

for the time-reversed process. Then, for $f, g \in C_0(E)$ and $h > 0$,

$$E(f(X_{t+h})g(X_t)) = E(g(X_t)E(f(X_{t+h})|X_t))$$
$$= E(P_h f(X_t) g(X_t))$$
$$= \int_E P_h f(x) g(x) \mu(dx).$$

We also have

$$E(f(X_{t+h})g(X_t)) = E(f(X_{t+h})E(g(X_t)|X_{t+h}))$$
$$= E(f(X_{t+h})\widetilde{P}_h g(X_{t+h}))$$
$$= \int_E f(x)\widetilde{P}_h g(x) \mu(dx).$$

We use that, due to the stationarity of the process $X$, $X_t$ and $X_{t+h}$ both have marginal distribution $\mu$. Combining these results gives

$$(3.2) \qquad \int_E P_h f(x) g(x) \mu(dx) = \int_E f(x) \widetilde{P}_h g(x) \mu(dx).$$

In the Appendix we calculate the operator $G^*$, defined by the equation

$$(3.3) \qquad \int_E Gf(x) g(x) \mu(dx) = \int_E f(y) G^* g(y) \mu(dy) \qquad \text{for all } f, g \in C_0(E).$$

It is shown there that

$$(3.4) \qquad G^* g(y) = \int_0^{T_1} g(\mathcal{L}_s y) \, ds + \frac{1}{\lambda} g(\mathcal{R} y) - \left(\frac{1}{\lambda} + T_1\right) g(y),$$

where in an analogous way as before we define $\mathcal{R} : E \to E$ as an exit to the right and $\mathcal{L}_s : E \to E$ as a new point at $s$ such that the point directly to the left of $s$ moves to the right.

We will use (3.4) several times. First of all, since $G^* 1 = 0$, it shows that $\mu$ is a stationary measure. Second, we see that for $g \in L^\infty(\mu)$

$$\|G^* g\|_\infty \leq 2\left(\frac{1}{\lambda} + T_1\right) \|g\|_\infty,$$

which proves that $G$ is in fact a continuous operator on $L^1(\mu)$, as well as a continuous operator on $C_0(E)$. Since $P_t = e^{tG}$, $P_t$ is also a continuous operator on $L^1(\mu)$. Therefore, (3.2) now shows that $\widetilde{P}_t = P_t^* = e^{tG^*}$, so in fact, using the same argument as before, $\widetilde{G} = G^*$. So the reversed process has the generator $G^*$.

Now we define a reflected Hammersley process $X^V$ as follows: we take the original stationary Hammersley process and reflect all the space–time paths with respect to the line segment $\{\frac{1}{2} T_1\} \times [0, T_2]$; call this a *vertical* reflection. So all points now move to the right and exit on the East side. One verifies that the generator for $X^V$ is given by $G^*$ in the same way we did it for the process $X$, and as $X^V$ also starts with a Poisson distribution of intensity $\lambda$, it has the same distribution as $\widetilde{X}$. Note that if one wishes to make a picture of the space–time paths of $\widetilde{X}$, one can take the original Hammersley process and reflect all the space–time paths with respect to the line-segment $[0, T_1] \times \{\frac{1}{2} T_2\}$, a *horizontal* reflection.

Since in $X^V$ all the jumps in $(0, T_1) \times (0, T_2)$ are made toward a point of a vertically reflected Poisson process, and in the process $\widetilde{X}$ all these jumps are made to the horizontally reflected $\beta$-points of the original Hammersley process, we have proved that the $\beta$-points are distributed according to a Poisson process with intensity 1. Furthermore, in the process $X^V$ paths exit on the East side according to a Poisson process with intensity $1/\lambda$, and this corresponds to $L_\lambda^{\text{in}}$, horizontally reflected. The process $L_\lambda^{\text{out}}$, also horizontally reflected, corresponds to the entries of $X^V$ at the $x$-axis, and is therefore Poisson with intensity $\lambda$. Finally, the independence of the three processes follows from the fact that this is true (by construction) for $X^V$. □

Theorem 3.1 allows us to show that a longest weakly North-East path from $(0,0)$ to $(t/\lambda^2, t)$ only spends a vanishing proportion of time on either the $x$- or $y$-axis. For the concept of longest weakly North-East path, see Definition 2.1.

COROLLARY 3.1. *Under the same conditions as Theorem* 3.1, *a longest weakly North-East path from* $(0,0)$ *to* $(t/\lambda^2, t)$ *spends a vanishing proportion of time on either the $x$- or $y$-axis, in the sense that the maximum distance from* $(0,0)$ *of the point where a longest weakly North-East path leaves the $x$- or $y$-axis, divided by $t$, tends to zero with probability* 1, *as* $t \to \infty$.

PROOF. Consider a longest weakly North-East path from $(0,0)$ to $(t/\lambda^2, t)$. Such a path can be associated with a path of a second-class particle from $(t/\lambda^2, t)$ to $(0,0)$ for the time-reversed process, running through the same $\alpha$-points as the longest weakly North-East path, but for which the roles of $\alpha$- and $\beta$-points are interchanged. This means that for the reversed process the associated path lies below or coincides with the path of the second-class particle that starts moving through the crossing of the upper edge $[0, t/\lambda^2] \times \{t\}$, closest to $(t/\lambda^2, t)$, moves down to the first $\alpha$-point on the path of the crossing, then moves to the left until it hits the path below the highest path crossing the rectangle $[0, t/\lambda^2] \times [0, t]$, then moves down again, and so on. Similarly this path lies above or coincides with the path of the second-class particle that starts moving to the left through the crossing of the right edge $\{t/\lambda^2\} \times [0, t]$, closest to $(t/\lambda^2, t)$, starts moving down when it hits the $\alpha$-point on the path of the crossing, moves to the left when it hits the next path, and so on.

According to Theorem 2.1 and Remark 2.1, now applied on the reversed process, the "$\beta$ waves" of the lower and upper path are asymptotically linear along the line through the origin with slope $\lambda^2$. This implies the statement of Corollary 3.1. □

REMARK 3.1. It is proved in Baik and Rains (2000) that $t^{-1/3}\{L_\lambda(t,t) - 2t\}$, where $L_\lambda(t,t)$ is the length of a longest North-East path from $(0,0)$ to $(t,t)$ in the stationary Hammersley process (as defined in Theorem 3.1, with $\lambda = 1$), converges in distribution to a distribution function $F_0$, which is related to, but different from the Tracy–Widom distribution function. This has the interesting consequence that the correlation between the number of points on the left edge and the number of crossings of the upper edge of the square $[0, t]^2$ tends to $-1$, as $t \to \infty$. Otherwise the variance of $L_\lambda(t,t)$ would be larger than $\eta t$, for some $\eta > 0$, instead of being of order $O(t^{2/3})$. We do not need their result in our argument, however. Baik and Rains (2000) use an analytical approach, applying the Deift–Zhou steepest descent method to an appropriate Riemann–Hilbert problem (after using a representation of the distribution function in terms of Toeplitz determinants). This approach is rather different from the approach taken here.

As noted in Baik and Rains (2000), the stationary process is a transition between two situations: if the intensities of the Poisson processes on the $x$-axis and $y$-axis are strictly smaller than 1, we get that $t^{-1/3}\{L_\lambda(t,t) - 2t\}$ converges in distribution to the Tracy–Widom distribution. On the other hand, if one of these intensities is bigger than 1 (but the intensities are not equal), we get convergence of $L_\lambda(t,t)$ to a normal distribution, with the usual $t^{-1/2}$ scaling (and a different centering constant).

REMARK 3.2. In Groeneboom (2001) a signed measure process $V_t$ was introduced, counting $\alpha$- and $\beta$-points contained in regions of the plane. The $V_t$-measure of a rectangle $[0,x] \times [0,y]$ is defined as the number of $\alpha$-points minus the number of $\beta$-points in the rectangle $[0,tx] \times [0,ty]$, divided by $t$. The $V_t$-process has the property that

$$V_t(S) \to V(S),$$

almost surely, for rectangles $S$ in the plane, where $V$ is a positive measure with density

$$(3.5) \qquad f_V(x,y) \stackrel{\text{def}}{=} \frac{\partial^2}{\partial x\, \partial y} V(x,y) = \frac{c}{4\sqrt{xy}}, \qquad x,y > 0.$$

Here we use the notation $V(x,y)$ to denote the $V$-measure of the rectangle $[0,x] \times [0,y]$. Likewise we write $V_t(x,y)$ for the $V_t$-measure of the rectangle $[0,x] \times [0,y]$.

The problem of proving part (ii) of Corollary 2.1 of the present paper was reduced to showing that

$$(3.6) \qquad \int_B \widetilde{V}_t(u,v)\, dV_t(u,v) \xrightarrow{\text{a.s.}} \int_B V(u,v)\, dV(u,v) = \tfrac{1}{4} c^2 xy,$$

where

$$\widetilde{V}_t(u,v) = \int_{[0,u] \times [0,v)} dV_t(u',v').$$

Although (3.6) indeed has to hold, the argument for it, given in Groeneboom (2001), is incomplete, and needs a result like Theorem 2.1 of the present paper to be completed. [The difficulty is caused by the locally unbounded variation of the measure $V_t$, as $t \to \infty$, which has to be treated carefully to explain why we need $\widetilde{V}_t$ as integrand in the integral in the left-hand side of (3.6) instead of, e.g., $V_t$, which leads to an integral that is asymptotically twice as large.] But since Theorem 2.1 allows us to prove both the local convergence to a Poisson process and convergence of $EL(t,t)/t$ to 2, we did not pursue the approach in Groeneboom (2001) any further in the present paper.

## APPENDIX

The purpose of this Appendix is to prove (3.4). Remember that

$$E = \bigsqcup_{n=0}^{\infty} E_n$$

where $E_0 = \{\varnothing\}$ and

$$E_n = \{(x_1, \ldots, x_n) : 0 \leq x_1 \leq \cdots \leq x_n \leq T_1\}.$$

A Poisson process of intensity $\lambda$ induces a probability measure $\mu$ on $E$. Denote by $\mu_n$ the restriction of $\mu$ to $E_n$, so $\mu_n(dx) = \lambda^n e^{-aT_1} dx$. The generator was given by

$$G : C_0(E) \to C_0(E) : Gf(x) = \int_0^{T_1} f(\mathcal{R}_t x) \, dt + \frac{1}{\lambda} f(\mathcal{L} x) - \left(\frac{1}{\lambda} + T_1\right) f(x).$$

Define $G_+ f = Gf + (1/\lambda + T_1)f$; we will calculate the dual of $G_+$. Let $f, g \in C_0(E)$:

$$\int_E G_+ f(x) g(x) \mu(dx)$$

$$= e^{-\lambda T_1} G_+ f(\varnothing) g(\varnothing) + \sum_{n=1}^{\infty} \int_{E_n} G_+ f(x) g(x) \mu_n(dx)$$

$$= e^{-\lambda T_1} \frac{1}{\lambda} f(\varnothing) g(\varnothing) + e^{-\lambda T_1} \int_0^{T_1} f(t) g(\varnothing) \, dt$$

$$+ e^{-\lambda T_1} \sum_{n=1}^{\infty} \left[\lambda^n \int_{E_n} \int_0^{T_1} f(\mathcal{R}_t x) g(x) \, dt \, dx + \lambda^{n-1} \int_{E_n} f(\mathcal{L} x) g(x) \, dx\right]$$

$$= e^{-\lambda T_1} \frac{1}{\lambda} f(\varnothing) g(\varnothing) + e^{-\lambda T_1} \int_0^{T_1} f(t) g(\varnothing) \, dt$$

$$+ e^{-\lambda T_1} \sum_{n=1}^{\infty} \sum_{i=1}^{n} \lambda^n \int_{\{x \in E_n, x_{i-1} < t \leq x_i\}} f(x_1, \ldots, x_{i-1}, t, x_{i+1}, \ldots, x_n)$$

$$\times g(x) \, dx \, dt$$

$$+ e^{-\lambda T_1} \sum_{n=1}^{\infty} \lambda^n \int_{\{x \in E_n, t > x_n\}} f(x_1, \ldots, x_n, t) g(x) \, dx \, dt$$

$$+ e^{-\lambda T_1} \sum_{n=1}^{\infty} \lambda^{n-1} \int_{E_n} f(x_2, \ldots, x_n) g(x) \, dx.$$

Now we make a change of variable for each term in such a way that we get $f(y)$ in each of the integrals:

$$\int_E G_+ f(x) g(x) \mu(dx)$$

$$= e^{-\lambda T_1} \frac{1}{\lambda} f(\varnothing) g(\varnothing) + e^{-\lambda T_1} \int_0^{T_1} f(y) g(\varnothing) \, dy$$

$$+ e^{-\lambda T_1} \sum_{n=1}^{\infty} \sum_{i=1}^{n} \lambda^n \int_{\{y \in E_n, y_i \leq s \leq y_{i+1}\}} f(y) g(y_1, \ldots, y_{i-1}, s,$$

$$y_{i+1}, \ldots, y_n) \, dy \, ds$$

$$+ e^{-\lambda T_1} \sum_{n=1}^{\infty} \lambda^n \int_{E_{n+1}} f(y) g(y_1, \ldots, y_n) \, dy$$

$$+ e^{-\lambda T_1} \sum_{n=1}^{\infty} \lambda^{n-1} \int_{\{y \in E_{n-1}, s \leq y_1\}} f(y) g(s, y_1, \ldots, y_{n-1}) \, dy \, ds$$

$$= \frac{1}{\lambda} f(\varnothing) g(\varnothing) \mu_0(E_0) + \frac{1}{\lambda} \int_{E_1} f(y) g(\varnothing) \mu_1(dy)$$

$$+ \sum_{n=1}^{\infty} \sum_{i=1}^{n} \int_{\{y \in E_n, y_i \leq s \leq y_{i+1}\}} f(y) g(y_1, \ldots, y_{i-1}, s,$$

$$y_{i+1}, \ldots, y_n) \mu_n(dy) \, ds$$

$$+ \sum_{n=0}^{\infty} \int_{\{y \in E_n, s \leq y_1\}} f(y) g(s, y_1, \ldots, y_n) \mu_n(dy) \, ds$$

$$+ \sum_{n=2}^{\infty} \frac{1}{\lambda} \int_{E_n} f(y) g(y_1, \ldots, y_{n-1}) \mu_n(dy)$$

$$= \sum_{n=0}^{\infty} \int_{E_n} f(y) \left( \int_0^{T_1} g(\mathcal{L}_s y) \, ds \right) \mu_n(dy) + \sum_{n=0}^{\infty} \frac{1}{\lambda} \int_{E_n} f(y) g(\mathcal{R} y) \mu_n(dy)$$

$$= \int_E f(y) \left( \int_0^{T_1} g(\mathcal{L}_s y) \, ds + \frac{1}{\lambda} g(\mathcal{R} y) \right) \mu(dy).$$

Here we define $\mathcal{R}$ as an exit to the right and $\mathcal{L}_s$ as a new point at $s$ such that the point directly to the left of $s$ moves to the right, that is,

$$\mathcal{R} : E \to E : \mathcal{R} x = \begin{cases} (x_1, \ldots, x_{n-1}), & \text{if } x \in E_n \ (n \geq 2), \\ \varnothing, & \text{if } x \in E_0 \sqcup E_1, \end{cases}$$

and for $0 < s < T_1$,

$$\mathcal{L}_s : E \to E : \mathcal{L}_s x = \begin{cases} (x_1, \ldots, x_{i-1}, s, x_{i+1}, \ldots, x_n), \\ \qquad \text{if } x_i \leq s < x_{i+1} \ (x \in E_n), \\ (s, x_1, \ldots, x_n), \qquad \text{if } s < x_1 \ (x \in E_n). \end{cases}$$

Since $G^* g = G^*_+ g - (1/\lambda + T_1) g$, we have shown that

$$G^* g(y) = \int_0^{T_1} g(\mathcal{L}_s y) \, ds + \frac{1}{\lambda} g(\mathcal{R}y) - \left(\frac{1}{\lambda} + T_1\right) g(y).$$

**Acknowledgments.** We are much indebted to Ronald Pyke for his comments and encouragement. We also want to thank Timo Seppäläinen for pointing out the connection of our result with the theory of second-class particles, which led to a simplification of the original proofs. Finally, we would like to thank an Associate Editor and referee for their helpful remarks.

## REFERENCES

ALDOUS, D. and DIACONIS, P. (1995). Hammersley's interacting particle process and longest increasing subsequences. *Probab. Theory Relatated Fields* **103** 199–213. MR1355056

ALDOUS, D. and DIACONIS, P. (1999). Longest increasing subsequences: From patience sorting to the Baik–Deift–Johansson theorem. *Bull. Amer. Math. Soc.* **36** 413–432. MR1694204

BAIK, J. and RAINS, E. (2000). Limiting distributions for a polynuclear growth model with external sources. *J. Statist. Phys.* **100** 523–541. MR1788477

BURKE, P. J. (1956). The output of a queueing system. *Oper. Res.* **4** 699–704. MR83416

FERRARI, P. A. (1992). Shocks in the Burgers equation and the asymmetric simple exclusion process. In *Automata Networks, Dynamical Systems and Statistical Physics* (E. Goles and S. Martinez, eds.) 25–64. Kluwer, Dordrecht. MR1263704

GROENEBOOM, P. (2001). Ulam's problem and Hammersley's process. *Ann. Probab.* **29** 683–690. MR1849174

GROENEBOOM, P. (2002). Hydrodynamical methods for analyzing longest increasing subsequences. *J. Comput. Appl. Math.* **142** 83–105. MR1910520

HAMMERSLEY, J. M. (1972). A few seedlings of research. *Proc. Sixth Berkeley Symp. Math. Statist. Probab.* **1** 345–394. Univ. California Press, Berkeley. MR405665

KINGMAN, J. F. C. (1973). Subadditive ergodic theory. *Ann. Probab.* **1** 883–909. MR356192

LIGGETT, T. M. (1999). *Stochastic Interacting Systems, Contact, Voter and Exclusion Processes.* Springer, New York. MR1717346

LOGAN, B. F. and SHEPP, L. A. (1977). A variational problem for random Young tableaux. *Adv. Math.* **26** 206–222. MR1417317

O'CONNELL, N. and YOR, M. (2002). A representation for non-colliding random walks. *Electron. Comm. Probab.* **7** 1–12. MR1887169

SEPPÄLÄINEN, T. (1996). A microscopic model for the Burgers equation and longest increasing subsequences. *Electron. J. Probab.* **1** 1–51. MR1386297

VERSHIK, A. M. and KEROV, S. V. (1977). Asymptotics of the Plancherel measure of the symmetric group and the limiting form of Young tableaux. *Soviet Math. Dokl.* **18** 527–531. (Translation of *Dokl. Acad. Nauk SSSR* **32** 1024–1027.) MR480398

26 E. CATOR AND P. GROENEBOOMDepartment of Applied Mathematics (DIAM)
Delft University of Technology
Mekelweg 4
2628 CD Delft
The Netherlands
e-mail: e.a.cator@ewi.tudelft.nl
e-mail: p.groeneboom@ewi.tudelft.nl